\documentclass[letterpaper, 10 pt, conference]{ieeeconf}  

\IEEEoverridecommandlockouts                              
\usepackage[utf8]{inputenc}
\usepackage[T1]{fontenc}

\usepackage{cite}
\usepackage{amsmath,amssymb,amsfonts}
\usepackage{algorithmic}
\usepackage{textcomp}
\usepackage{color}
\usepackage{bm}
\usepackage{graphicx}
\usepackage{enumerate}
\usepackage{multicol}

\newtheorem{lemma}{Lemma}
\newtheorem{assumption}{Assumption}
\newtheorem{remark}{Remark}
\newtheorem{theorem}{Theorem}
\newtheorem{example}{Example}

\title{
\LARGE \bf
Robust Accelerated Dynamics for Subnetwork Bilinear Zero-Sum Games with Distributed Restarting
}

\author{Weijian~Li, Lacra Pavel, and Andreas A. Malikopoulos
\thanks{
This research was supported in part by NSF under Grants CNS-2401007, CMMI-2348381, IIS-2415478, and in part by MathWorks.}
\thanks{W.~Li and A.A. Malikopoulos are with the Department of Civil and Environmental Engineering, Cornell University, Ithaca, NY, USA (emails: \texttt{wl779@cornell.edu},
\texttt{amaliko@cornell.edu}).}
\thanks{L. Pavel is with the Department of Electrical and Computer Engineering, University of Toronto, Toronto, ON, Canada
(email: \texttt{pavel@control.utoronto.ca}).
}
}

\begin{document}

\maketitle

\begin{abstract}
In this paper, we investigate distributed Nash equilibrium seeking for a class of two-subnetwork zero-sum games characterized by bilinear coupling. We present a distributed primal-dual accelerated mirror-descent algorithm with convergence guarantees. However, we demonstrate that this time-varying algorithm is not robust, as it fails to converge under arbitrarily small disturbances. To address this limitation, we introduce a distributed accelerated algorithm 
that incorporates a coordinated restarting mechanism.
We model this new algorithm as a hybrid dynamical system and establish its structural robustness.
\end{abstract}


\section{INTRODUCTION}
In recent years, there has been a notable surge of research interest in distributed Nash equilibrium (NE) seeking for noncooperative games,
driven by their broad applications in various areas such as sensor networks, smart grids, and social networks \cite{pavel2022distributed, ye2023distributed, chremos2024mechanism}.
In particular, the two-subnetwork (bilinear) zero-sum game is one of the most important classes of noncooperative games.
To be specific, agents in one subnetwork aim to minimize a global payoff function using their local information, whereas agents in the other subnetwork seek to maximize the same function. Multi-agent decision-making problems, such as coordination of autonomous vehicles and power allocation in the presence of adversaries, can be effectively modeled by subnetwork zero-sum games \cite{Malikopoulos2021,Dave2024TAC, chalaki2020TITS, gharesifard2013distributed}.
Therefore, it is worthwhile to design efficient distributed algorithms for NE seeking in such games.

First-order primal-dual methods have been extensively investigated for subnetwork zero-sum games since they are powerful in handling constraints, suitable for large-scale problems, and easily implemented in a distributed manner.
For instance, distributed continuous-time primal-dual algorithms were proposed in \cite{gharesifard2013distributed}, while a discrete-time subgradient-based algorithm was designed in \cite{lou2015nash} under switching topologies. For nonsmooth payoff functions and bounded constraints, a distributed algorithm was explored in \cite{yang2018cooperative}.
However, standard primal-dual methods can only achieve convergence with a rate of $\mathcal O(1/t)$ for general convex cost functions \cite{zeng2018distributed}.
Inspired by momentum-based acceleration protocols such as the heavy-ball algorithm and Nesterov's scheme \cite{nesterov2003introductory},  efforts have been made to develop accelerated algorithms under a primal-dual framework.
An accelerated primal-dual algorithm with adaptive parameters was explored for strongly convex cost functions in \cite{xu2017accelerated}, while an accelerated protocol with an optimal linear rate was designed in \cite{kovalev2022accelerated}.
Continuous-time primal-dual Nesterov accelerated dynamics were developed in \cite{zeng2022dynamical, zeng2023distributed}.
In \cite{li2024primal}, a primal-dual accelerated mirror-descent algorithm was proposed in order to handle constraint sets.

Momentum-based algorithms may achieve faster convergence rates, but at the expense of sacrificing the structural robustness properties \cite{o2015adaptive, devolder2014first}. 
The authors of \cite{poveda2019inducing} showed that continuous-time Nesterov's accelerated dynamics were not uniformly globally asymptotically stable, and could not be convergent in the presence of disturbances.
To ensure convergence and robustness, a hybrid dynamical system (HDS) with a restarting mechanism was proposed in \cite{poveda2019inducing}.
Following that, hybrid zero-order algorithms were designed in \cite{poveda2021robust}.
The ideas were generalized to distributed optimization in \cite{ochoa2020robust} and NE seeking in \cite{ochoa2023momentum}.
However, the acceleration protocol in \cite{poveda2019inducing, poveda2021robust, ochoa2020robust, ochoa2023momentum}  was limited to Nesterov's scheme, and could only solve unconstrained optimization/game problems.

This paper explores a distributed accelerated algorithm for subnetwork zero-sum bilinear games.
Our main contributions are twofold. 
a) The formulation is a generalization of the problems in \cite{xu2017accelerated, kovalev2022accelerated, zeng2022dynamical} by allowing constraint sets and strictly convex cost functions. We revisit the convergence of a distributed primal-dual accelerated mirror-descent algorithm in \cite{li2024primal}, and find that it fails to converge in the presence of arbitrarily small disturbances. 
b) By incorporating a coordinated restarting mechanism into the time-varying mirror-descent method, we develop a novel distributed accelerated algorithm, which covers those in \cite{ochoa2020robust, ochoa2023momentum} as special cases.
We establish that the algorithm reaches an NE with structural robustness properties.


\section{MATHEMATICAL PRELIMINARIES}

We now introduce the notation and preliminary results that will be used throughout the paper.

\emph{Notation:}
Let $\mathbb R_{\ge 0}$ and $\mathbb Z_{\ge 0}$ be the set of nonnegative real numbers and nonnegative integers.
Denote by $0_p$ ($1_p$) the $p$-dimensional column vector with all entries of $0$ ($1$), and by $I_p$ the $p$-by-$p$  identity matrix. 
We use $\bm 0$ to denote vectors/matrices of zeros with appropriate dimensions when there is no confusion.
Denote by $(\cdot)^\top$, $\otimes$, and $\Vert\cdot\Vert$ the transpose, the Kronecker product, and the Euclidean norm, respectively.
For $x_i \in \mathbb R^{m_i}$, ${\rm col}\{x_1, \dots, x_n\}=[x_1^\top, \dots, x_n^\top]^\top$.
Let $\rho \mathbb B$ be the closed ball of appropriate dimension in the Euclidean norm of radius $\rho$.
Define $|x|_{\mathcal A} = \min_{y \in \mathcal A} \Vert x - y\Vert$ for a compact set $\mathcal A$.
Given a differential function $f: \mathbb R^p \rightarrow \mathbb R$, $\nabla f(x)$ is the gradient $f$ at $x$.
A continuous function $\alpha: \mathbb R_{\ge 0} \rightarrow \mathbb R_{\ge 0}$ is of class $\mathcal{K}$ if $\alpha(0) = 0$ and strictly increasing, and it is of class $\mathcal{K}_\infty$ if it is of class $\mathcal{K}$, and $\alpha(r) \rightarrow \infty$ as $r \rightarrow \infty$.
A continuous function $\beta: \mathbb R_{\ge 0} \times\mathbb R_{\ge 0} \rightarrow \mathbb R_{\ge 0}$ is of class $\mathcal{KL}$ if $\beta(r, s)$ belongs to class $\mathcal K$ with respect to $r$ for each fixed $s$,  $\beta(r, s)$ is decreasing with respect to $s$ for each fixed $r$, and $\beta(r, s) \rightarrow 0$ as $s \rightarrow \infty$.

\emph{Convex Analysis:} 
Let $\Omega \subset \mathbb R^p$ be a closed convex set, and 
$f: \Omega \rightarrow \mathbb R$ be a continuously differentiable function. 
Then $f$ is convex if 
$f(y) \ge f(x) + \langle y - x, \nabla f(x)\rangle, \forall x, y \in \Omega$,
strictly convex if the strict inequality holds for all $x \not= y$,
and strongly convex if
$f(y) \ge f(x) + \langle y - x, \nabla f(x)\rangle + \frac \mu 2 \Vert x - y\Vert^2, \forall x, y \in \Omega$ for some $\mu > 0$.
Suppose that $f$ is convex. Then its convex conjugate is 
$f^*(u) = \sup_{x \in \Omega}\{x^\top u - f(x)\}$, and its Bregman divergence is $D_f(x, y) = f(x) - f(y) - \langle x - y, \nabla f(y)\rangle$.
In light of \cite{diakonikolas2019approximate}, the following result holds.
\begin{lemma}
Let $\Omega$ be a closed convex set and $f: \Omega \rightarrow \mathbb R$ be a differentiable and strongly convex function. Then $f^*$ is differentiable and convex such that $\nabla f^*(u) = {\rm argmin}_{x \in \Omega}\{-x^\top u + f(x)\}$.
\end{lemma}

\emph{Hybrid Dynamical System:}
Referring to \cite{goedel2012hybrid}, we model a HDS as
\begin{equation}
\label{HDS}
x \in \mathcal C, ~\dot x = F(x),~{\rm and}~ 
x \in \mathcal D, ~x^+ \in G(x),
\end{equation}
where $x \in \mathbb R^n$ is the state, 
$\mathcal C \subset \mathbb R^n$ is the flow set, 
$F: \mathbb R^n \rightarrow \mathbb R^n$ is the flow map,
$\mathcal D \subset \mathbb R^n$ is the jump set, and
$G: \mathbb R^n \rightrightarrows \mathbb R^n$ is the set-valued jump map.
We use $\mathcal H = \{\mathcal C, F, \mathcal D, G\}$ to denote the HDS (\ref{HDS}).
Solutions to $\mathcal H$ are defined on hybrid time domains, i.e., they are parameterized by a continuous-time index $t \in \mathbb R_{\ge 0}$ and a discrete-time index $j \in \mathbb Z_{\ge 0}$.
A compact set $\mathcal A \subset C \cup D$ is uniformly globally asymptotically stable (UGAS) if 
$\exists \beta \in \mathcal {KL}$ such that every solution $x$ satisfies
$|x(t, j)|_{\mathcal A} \le \beta(|x(0, 0)|_{\mathcal A}, t + j), \forall (t, j) \in {\rm dom}(x).$

Consider a $\epsilon$-parameterized HDS $\mathcal H_{\epsilon}$ as
$x \in \mathcal C_\epsilon, \dot x = F_\epsilon(x)$, and 
$x \in \mathcal D_\epsilon, \dot x \in G_\epsilon(x)$, 
where $\epsilon > 0$. The compact set $\mathcal A$ is semi-globally practically asymptotically stable (SGPAS) as
$\epsilon \to 0^+$ with $\beta \in \mathcal {KL}$ if for each pair $\delta > \mu > 0$, there exists  $\epsilon^*$ such that for all $\epsilon \in (0, \epsilon^*)$, every solution $x$ of $\mathcal H_{\epsilon}$ with $|x(0, 0)|_{\mathcal A} \le \delta$ satisfies
$|x(t, j)|_{\mathcal A} \le \beta(|x(0, 0)|_{\mathcal A}, t + j) + {\mu}.$

Consider a perturbed HDS given by
\begin{equation}
\begin{aligned}
\label{HDS:D}
&x + e \in \mathcal C, &&\dot x = F(x + e) + e, \\
&x + e \in \mathcal D, &&x^+ \in G(x + e) + e,
\end{aligned}
\end{equation}
where $e$ is measurable, and $\sup_{(t, j) \in {\rm dom}(e)} |e(t, j)| \le \epsilon$ for some $\epsilon > 0$. System (\ref{HDS:D}) is UGAS with a structurally robust property (R-UGAS) if it is UGAS when $\epsilon = 0$, and meanwhile, it is SGPAS as $\epsilon \to 0^+$.

\emph{Graph Theory:}
A multi-agent network can be modeled by an undirected graph $\mathcal G(\mathcal V, \mathcal E)$, where $\mathcal V = \{1, \dots, n\}$ is the node set and $\mathcal E \subset \mathcal V \times \mathcal V$ is the edge set. The graph $\mathcal G$ is associated with a weight matrix $A = [a_{ij}] \in \mathbb R^{n \times n}$ 
such that $a_{ij} = a_{ji} > 0$ if $(i, j) \in \mathcal E$, and $a_{ij} = 0$ otherwise. Its Laplacian matrix $\mathcal L$ is defined as $\mathcal L = D - A$, where $D = {\rm diag}\{d_i\}$, and $d_i = \sum_{j \in \mathcal V} a_{ij}$. The graph $\mathcal G$ is connected if there exists a path between any pair of distinct nodes.

\section{PROBLEM STATEMENT}

Consider a network $\mathcal G$ consisting of two undirected and connected subnetworks  $\mathcal G_1(\mathcal V_1, \mathcal E_1)$ and $\mathcal G_2(\mathcal V_2, \mathcal E_2)$ with 
$\mathcal V_1 = \{1, \dots, n_1\}$ and 
$\mathcal V_2 = \{1, \dots, n_2\}$.
The strategic variables of subnetworks $\mathcal G_1$ and $\mathcal G_2$ are
$\bm x = {\rm col}\{x_1, \dots, x_{n_1}\}$ and
$\bm y = {\rm col}\{y_1, \dots, y_{n_2}\}$,
and the feasible sets are 
$\Omega_1 = \{\bm x \in X: x_1 = \dots = x_{n_1}\}$ and 
$\Omega_2 = \{\bm y \in Y: y_1 = \dots = y_{n_2}\}$,
where $X = X_1 \times \dots \times X_{n_1}$,
$Y = Y_1 \times \dots \times Y_{n_2}$, and moreover,
$X_i \subset \mathbb R^{p_1}$ and
$Y_j \subset \mathbb R^{p_2}$ are closed convex sets for all 
$i \in \mathcal V_1$ and
$j \in \mathcal V_2$.
The payoff function of $\mathcal G_1$ ($\mathcal G_2$) is
$U$ ($-U$), where $U: X \times Y \rightarrow \mathbb R$ is 
$$
U(\bm x, \bm y) := 
\tilde f(\bm x) + \bm x^\top B \bm y - \tilde g(\bm y),
$$
with $\tilde f(\bm x) = \sum_{i \in \mathcal V_1} f_i(x_i)$,
$\tilde g(\bm y) = \sum_{j \in \mathcal V_2} g_j(y_j)$ and
$\bm x^\top B \bm y = \sum_{i \in \mathcal V_1}\sum_{j \in \mathcal V_2} x_i^\top H_{ij} y_j$.
Let $\mathcal E_{\bowtie}$ be the set of connected edges between $\mathcal G_1$ and $\mathcal G_2$, i.e., agent $i$ of $\mathcal G_1$ and agent $j$ of $\mathcal G_2$ can observe each other's decision variable if $(i, j) \in \mathcal E_{\bowtie}$.
Suppose $H_{ij} \in \mathbb R^{p_1 \times p_2}$ is a non-zero matrix only if 
$(i, j) \in \mathcal E_{\bowtie}$.
The agents in $\mathcal G_1$ compute $\bm x$ to minimize $U$, while the agents in $\mathcal G_2$ find $\bm y$ to maximize $U$.
However, agent $i$ in $\mathcal G_1$ only knows $f_i$, $X_i$ and $H_{ij}$, and agent $j$ in $\mathcal G_2$ knows $g_j$, $Y_j$ and $H_{ij}$, where $(i, j) \in \mathcal E_{\bowtie}$.
Seeking an NE of the game can be cast into
\begin{equation}
\begin{aligned}
\label{formulation}
\min_{\bm x \in X} \max_{\bm y \in Y}~ U(\bm x, \bm y),
~~{\rm s.t.}~ \bm L_1 \bm x =\bm 0, ~ 
\bm L_2 \bm y = \bm 0,
\end{aligned}
\end{equation}
where $\bm L_1 = \mathcal L_1 \otimes I_{p_1}$,
$\bm L_2 = \mathcal L_2 \otimes I_{p_2}$, and
$\mathcal L_1$, $\mathcal L_2$ are the Laplacian matrices of $\mathcal G_1$ and $\mathcal G_2$.

To ensure the well-posedness and solvability of the game, we impose the following assumption \cite{gharesifard2013distributed, zeng2022dynamical, gao2022continuous}.
\begin{assumption}
\label{ass:convex}
The functions $f_i$ and $g_j$ are continuously differentiable and strictly convex on $X_i$ and $Y_j$, respectively.
Their gradients are Lipschitz continuous over $X_i$ and $Y_j$. 
The Slater’s constraint qualification holds.
Moreover, the NE lies in the relative interior of $X \times Y$.
\end{assumption}

\begin{remark}
The formulation \eqref{formulation} covers those in \cite{zeng2023distributed, kovalev2022accelerated} as special cases by allowing strictly convex cost functions and constraint sets.
In practice, \eqref{formulation} appears in distributed optimization with parameter uncertainties and adversarial resource allocation of multiple communication channels \cite{gharesifard2013distributed}.
\end{remark}

Following \cite{gharesifard2013distributed}, we equivalently reformulate (\ref{formulation}) as
\begin{equation}
\begin{aligned}
\label{lagrangian}
\min_{\bm x, \bm \mu} ~
\max_{\bm y, \bm \lambda}~
S(\bm x, \bm \lambda, \bm y, \bm \mu)
: = U(\bm x, \bm y) + \bm \lambda^\top \bm L_1 \bm x& \\
- \bm \mu^\top \bm L_2 \bm y + \frac 12 \bm x^\top \bm L_1 \bm x - \frac 12 \bm y^\top \bm L_2 \bm y,
\end{aligned}
\end{equation}
where $\bm \lambda ={\rm col}\{\lambda_1, \dots, \lambda_{n_1}\} \in \mathbb R^{p_1 n_1}$ and $\bm \mu = {\rm col}\{\mu_1, \dots, \mu_{n_2}\} \in \mathbb R^{p_2 n_2}$ are the Lagrangian multipliers.

Similar to the method in \cite{li2024primal}, we introduce auxiliary variables 
$\bm u = {\rm col}\{u_1, \dots, u_{n_1}\}$,
$\bm \gamma = {\rm col}\{\gamma_1, \dots, \gamma_{n_1}\}$,
$\bm v = {\rm col}\{v_1, \dots,v_{n_2}\}$ and
$\bm \nu = {\rm col}\{\nu_1, \dots, \nu_{n_2}\}$,
and generating functions
$\bm\psi(\bm x) = {\rm col}\{\psi_1(x_1), \dots, \psi_{n_1}(x_{n_1})\}$ and
$\bm\phi(\bm y) = {\rm col}\{\phi_1(y_1), \dots, \phi_{n_2}(y_{n_2})\}$.
Then we propose a distributed primal-dual accelerated mirror-descent algorithm as
\begin{equation}
\begin{aligned}
\label{alg}
\dot x_i &= \frac rt \big[\nabla \psi_i^*(u_i) - x_i \big], ~
\dot \lambda_i = \frac rt \big[\gamma_i - \lambda_i \big], \\
\dot u_i &= \frac tr \big[ - \nabla_{x_i} S\big(\bm x, \bm \lambda \! + \! \frac t r \dot{\bm \lambda}, \bm y + \frac t r \dot{\bm y}, \bm \mu \big) \big], \\
\dot \gamma_i &= \frac tr \big[\nabla_{\lambda_i} S \big(\bm x + \frac tr \dot{\bm x}, \bm \lambda, \bm y, \bm \mu \big)\big], \\
\dot y_j &= \frac rt \big[\nabla \phi_j^*(v_j) - y_j \big], ~
\dot \mu_j = \frac rt \big[\nu_j - \mu_j \big], \\
\dot v_j &= \frac tr \big[\nabla_{y_j}  S \big(\bm x  +  \frac t r \dot{\bm x},  \bm \lambda, \bm y, \bm \mu  +  \frac t r \dot{\bm \mu} \big) \big], \\
\dot \nu_j &= \frac tr \big[ -\nabla_{\mu_i} S \big(\bm x, \bm \lambda, \bm y  + \frac tr \dot{\bm y}, \bm \mu \big) \big],
\end{aligned}	
\end{equation}
where $t_0 > 0$,
$\psi_i^*$, $\phi_j^*$ are the convex conjugate functions of 
$\psi_i$ and $\phi_j$,
$\nabla \bm \psi^*(\bm u) = {\rm col}\{\nabla \psi_1^*(u_1), \dots, \nabla \psi^*_{n_1}(u_{n_1})\}$,
$\nabla \bm \phi^*(\bm v) = {\rm col}\{\nabla \phi_1^*(v_1), \dots, \nabla \phi^*_{n_1}(v_{n_2})\}$,
$\bm u(t_0) = \bm u_0 \in X$,
$\bm x(t_0) = \bm x_{0} = \nabla \bm\psi^*(\bm u_0)$,
$\bm \gamma(t_0) = \bm \gamma_0 \in \mathbb R^{p_1n_1}$,
$\bm \lambda(t_0) = \bm \lambda_0 = \bm \gamma_0$, 
$\bm v(t_0) = \bm v_0 \in Y$,
$\bm y(t_0) = \bm y_0 = \nabla \bm \phi^*(\bm v_0)$,
$\bm \nu(t_0) = \bm \nu_0 \in \mathbb R^{p_2n_2}$, and
$\bm \mu(t_0) = \bm\mu_0 = \bm \nu_0$.

We make the following assumption for $\bm\psi$ and $\bm \phi$ \cite{krichene2015accelerated, gao2022continuous}.
\begin{assumption}
\label{ass:genlip}
The generating functions $\psi_i: X_i \rightarrow \mathbb R$ and $\phi_j: Y_j \rightarrow \mathbb R$ are continuously differentiable with Lipschitz continuous gradients, and they are strongly convex.
Besides, $\nabla \psi_i^*$ and $\nabla \phi_j^*$ are Lipschitz continuous over $X_i$ and $Y_j$.
\end{assumption}

\begin{remark}
The dynamics (\ref{alg}) is a generalization of the algorithm in \cite{krichene2015accelerated} under a primal-dual framework. The variables $\bm u$ and $\bm \nu$ descend in the negative gradient direction of $S$, while $\bm v$ and $\bm \gamma$ ascend in the other direction. 
The derivative feedback, introduced as damping terms, is beneficial for its convergence. 
The time-varying coefficients $t/r$ and $r/t$ enable (\ref{alg}) to achieve acceleration.
In the absence of $X$ and $Y$, (\ref{alg}) degenerates into the  algorithm in \cite{zeng2023distributed, zeng2022dynamical}. 
\end{remark}

\begin{remark}
The generating functions are introduced to handle constraints and to analyze the convergence from the dual space \cite{ diakonikolas2019approximate}.
In practice, their selection depends on the constraint sets \cite{gao2022continuous, krichene2015accelerated}.
For instance, if $X_i$ is a closed convex set, we take $\psi_i(x_i) = \frac 12 \Vert x_i\Vert^2, x_i \in X_i$, and then, $\nabla \psi_i^*(u_i) = {\rm proj}_{X_i}(u_i)$, where ${\rm proj}_{X_i}(\cdot)$ is the Euclidean projection operator. 
\end{remark}

\begin{remark}
From Assumptions \ref{ass:convex} and \ref{ass:genlip}, it follows that there is a unique NE for (\ref{formulation}), and that there is a unique solution to (\ref{alg}). 
In light of \cite[Lemma 4]{li2024primal}, $x(t) \in X$ and $y(t) \in Y$ for all $t \ge 0$ in (\ref{alg}).
Let $(\bm u^\star, \bm x^\star, \bm \gamma^\star, \bm \lambda^\star, \bm v^\star, \bm y^\star, \bm \nu^\star, \bm \mu^\star)$ be an equilibrium point of (\ref{alg}).
By the Karush–Kuhn–Tucker (KKT) optimality conditions,
$(\bm x^\star, \bm \lambda^\star, \bm y^\star, \bm \mu^\star)$ is a saddle point to $S$, and $(\bm x^\star, \bm y^\star)$ is the NE.
Note that (\ref{alg}) is simpler than (15) in \cite{li2024primal} since extra terms in the evolution of $\bm u, \bm \gamma, \bm v$ and $\bm \mu$ are removed, but its states are still convergent because $(\bm x^\star, \bm y^\star)$ lies in the relative interior of $X \times Y$.
\end{remark}

The next theorem addresses the convergence of (\ref{alg}).
\begin{theorem}
\label{thm:alg}
Consider dynamics (\ref{alg}) with $r \ge 2$.
Let $\big(\bm x(t), \bm \lambda(t), \bm y(t), \bm \mu(t)\big)$ be a trajectory of (\ref{alg}), and $(\bm x^\star, \bm \lambda^\star, \bm y^\star, \bm \mu^\star)$ be a saddle point of $S$. Under Assumptions \ref{ass:convex} and \ref{ass:genlip}, the duality gap, defined by $S \big(\bm x(t), \bm \lambda^\star, \bm y^\star, \bm \mu (t)\big) - S(\bm x^\star, \bm \lambda(t), \bm y(t), \bm \mu^\star)$, converges to $0$ with a rate of $\mathcal O(1/t^2)$. Moreover, the trajectory of $\big(\bm x(t), \bm y(t)\big)$ approaches $(\bm x^\star, \bm y^\star)$.
\end{theorem}

\emph{Proof:} 
Define
\begin{equation}
\begin{aligned}
\label{pf:lya0}
V =& \frac {t^2}{r}\big(S(\bm x, \bm \lambda^\star, \bm y^\star, \bm \mu) - S(\bm x^\star, \bm \lambda, \bm y, \bm \mu^\star) \big) + r D_{\psi_{u}^*}(\bm u, \bm u^\star) \\
& + r D_{\psi_{\gamma}^*}(\bm \gamma, \bm \gamma^\star) 
+ r D_{\phi_{v}^*}(\bm v, \bm v^\star)
+ r D_{\phi_{\nu}^*}(\bm \nu, \bm \nu^\star),
\end{aligned}
\end{equation}
where $S$ is defined by (\ref{lagrangian}),
$(\bm u^\star, \bm x^\star, \bm \gamma^\star, \bm \lambda^\star, \bm v^\star, \bm y^\star, \bm \nu^\star, \bm \mu^\star)$ is an equilibrium point of (\ref{alg}),
$D_{\psi_{u}^*}(\bm u, \bm u^\star) = \sum_{i \in \mathcal V_1} D_{\psi_{i}^*}(u_i, u_i^\star)$,
$D_{\psi_{\gamma}^*}(\bm \gamma, \bm \gamma^\star) = \sum_{i \in \mathcal V_1} \frac 12 \Vert \gamma_i - \gamma_i^* \Vert^2$,
$D_{\phi_{v}^*}(\bm v, \bm v^\star)  = \sum_{j \in \mathcal V_2} 
D_{ \phi_j^*}(v_j, v_j^\star)$, and
$D_{\phi_{\nu}^*}(\bm \nu, \bm \nu^\star)  = \sum_{j \in \mathcal V_2} \Vert \mu_j - \mu_j^* \Vert^2$.

Plugging \eqref{alg} into \eqref{pf:lya0}, we obtain
\begin{equation*}
\begin{aligned}
\dot V &= \frac {2t}{r}\big(S(\bm x, \bm \lambda^\star, \bm y^\star, \bm \mu) - S(\bm x^\star, \bm \lambda, \bm y, \bm \mu^\star)\big) \\
&- t \big(\langle \bm x - \bm x^*, \nabla \tilde f(\bm x) 
+ \bm L_1 \bm x \rangle + {\bm\lambda}^{*, \top} \bm L_1 \bm x + \bm x^\top B \bm y^* \\
&~~+\langle \bm y - \bm y^*, \nabla \tilde g(\bm y)
+ \bm L_2 \bm y \rangle + {\bm\mu}^{*, \top} \bm L_2 \bm y - \bm x^{*, T} B \bm y\big).
\end{aligned}	
\end{equation*}

Since $S$ is convex with respect to $(\bm x, \bm \mu)$, and concave with respect to $(\bm y, \bm \lambda)$, we obtain
\begin{equation}
\label{ine:lya}
\dot V \le (\frac {2t}{r} - t)\big(S(\bm x, \bm \lambda^\star, \bm y^\star, \bm \mu) - S(\bm x^\star, \bm \lambda, \bm y, \bm \mu^\star)\big).
\end{equation}
Clearly, $\dot V \le 0$ if $r \ge 2$.
Let $m_0 = S(\bm x_0, \bm \lambda_0, \bm y_0, \bm \mu_0)$. Then
\begin{equation*}
\begin{aligned}
\label{pf:ine}
S(\bm x, \bm \lambda^\star, \bm y^\star, \bm \mu) - S(\bm x^\star, \bm \lambda, \bm y, \bm \mu^\star) \le \frac {r}{t^2} V \le \frac {r}{t^2} m_0,
\end{aligned}
\end{equation*}
i.e.,
$S(\bm x, \bm \lambda^\star, \bm y^\star, \bm \mu) - S(\bm x^\star, \bm \lambda, \bm y, \bm \mu^\star)$
converges to $0$ with a rate of $\mathcal O(1/t^2)$.

Due to the strict convexity of $\tilde f$ and $\tilde g$, $S(\bm x, \bm \lambda^\star, \bm y^\star, \bm \mu) = S(\bm x^\star, \bm \lambda, \bm y, \bm \mu^\star)$ only if $(\bm x, \bm y) = (\bm x^*, \bm y^*)$. 
Thus,
$(\bm x, \bm y)$ approaches $(\bm x^*, \bm y^*)$, and
the proof is completed.
$\hfill\square$

Theorem \ref{thm:alg} indicates that (\ref{alg}) reaches the NE with a rate of $\mathcal O(1/t^2)$, which is faster than the standard primal-dual method \cite{zeng2018distributed}.
Note that the requirement of $r \ge 2$ is typical for the convergence of accelerated algorithms \cite{su2016differential, krichene2015accelerated}.
However, as shown in \cite{poveda2019inducing, devolder2014first}, momentum-based algorithms may lack robustness. The following example is provided for illustration.
\begin{example}
Consider problem (\ref{formulation}) with
$n_1 = n_2 = 4$ and $p_1 = p_2 = 2$. 
Let 
$f_i(x_i) = \log[\exp(x_{i1}-0.1i)+\exp(x_{i2} - 0.2i)]$, 
$g_j(y_j) = \exp(y_{j1} - j) - y_{j1} + (y_{j2} - 0.3j)^2$,
and $B = I_{8}$. 
Take $X_i = \{x_i \in \mathbb R^2~|~l_i \le x_i \le u_i\}$ and $Y_i = X_i$, where $l_j, u_j \in \mathbb R^2$ are random vectors between $[-2, 0]$ and $[1, 3]$.
Fig. \ref{fig:Ex1} shows the trajectories of $U(\bm x, \bm y)$
under dynamics (\ref{alg}). Clearly, the dynamics reaches the NE as proved in Theorem \ref{thm:alg}.
However, when even a small disturbance $e = 0.001$ is added as in (\ref{HDS:D}), the convergence is lost.

\begin{figure}[htp]
\centering
\includegraphics[scale=0.4]{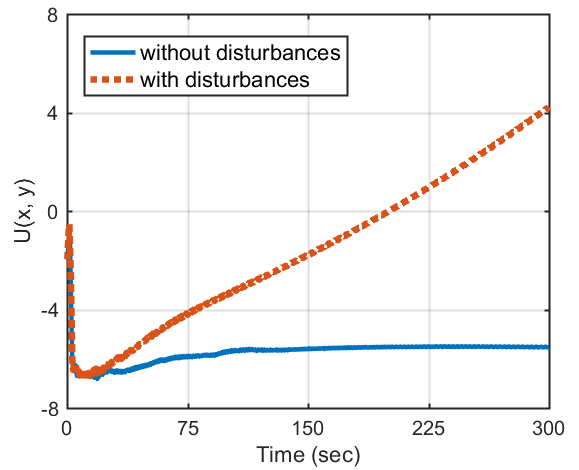}
\caption{Trajectories of $U(\bm x, \bm y)$ under dynamics (\ref{alg}).}
\label{fig:Ex1}
\end{figure}
\end{example}

The observation motivates us to develop an accelerated algorithm that enjoys both convergence and structural robustness.
We note that an approach as in \cite{romano2019dynamic} cannot be used, as an exact model of the disturbance in (\ref{HDS:D}) is unknown.

\section{DISTRIBUTED ACCELERATED ALGORITHM WITH COORDINATED RESTARTING}

This section introduces a robust accelerated algorithm for (\ref{formulation}) by incorporating a restarting technique into (\ref{alg}). 
This approach is inspired by \cite{poveda2019inducing, poveda2021robust, ochoa2020robust,ochoa2023momentum}, which showed that Nesterov's scheme, used with restarting techniques, may achieve R-UGAS.

\subsection{Hybrid Dynamics with Restarting}

Define 
$\bm\xi^1 = {\rm col}\{\bm x, \bm \lambda\}$,
$\bm\zeta^1 = {\rm col}\{\bm u, \bm \gamma\}$,
$\nabla \bm\pi^{1, *}(\bm\zeta^1) = {\rm col}\{\nabla \bm\psi^*(\bm u), \bm \gamma\}$,
$\nabla \bm h^1(\bm \xi^1) = {\rm col}\{\nabla \tilde f(\bm x) + \bm L_1 \bm x, \bm 0\}$,
$\tilde L^1 = [\bm 0, \bm L_1; - \bm L_1, \bm 0]$,
$\tilde B^1 = [B, \bm 0; \bm 0, \bm 0]$,
$\bm \xi^2 = {\rm col}\{\bm y, \bm \mu\}$,
$\bm \zeta^2 = {\rm col}\{\bm v, \bm \nu\}$,
$\nabla \bm\pi^{2, *}(\bm \zeta^2) = {\rm col}\{\nabla \bm \phi^*(\bm v), \bm\nu\}$,
$\nabla \bm h^2(\bm \xi^2) = {\rm col}\{\nabla \tilde g(\bm y) + \bm L_2 \bm y, \bm 0\}$,
$\tilde L^2 = [\bm 0, \bm L_2; - \bm L_2, \bm 0]$, and
$\tilde B^2 = [-B^\top, \bm 0; \bm 0, \bm 0]$.
Let
$\mathbb D_{\xi}^1 = X \times \mathbb R^{p_1 n_1}$,
$\mathbb D_{\xi}^2 = Y \times \mathbb R^{p_2 n_2}$ and
$\mathbb D_{\zeta}^l = \mathbb R^{2 p_l n_l}$.
Then dynamics (\ref{alg}) reads as
\begin{equation}
\begin{aligned}
\label{alg:asim}
\dot {\bm \xi}^l &= \frac {r}{t} \big[\nabla \bm\pi^{l, *} (\bm\zeta^l) - \bm\xi^l \big], \\
\dot {\bm \zeta}^l &= \frac tr \big[\!-\!\nabla \bm h^l(\bm\xi^l) \!-\! \tilde L^l (\bm\xi^l + \frac tr \dot {\bm\xi}^l) - \tilde B^l(\bm\xi^{(3 - l)} + \frac tr \dot {\bm\xi}^{(3 - l)}) \big],
\end{aligned}
\end{equation}
where $l \in \{1, 2\}$,
$\bm \zeta^l(t_0) = \bm \zeta^l_0 \in \mathbb D^l_{\bm \xi}$ and 
$\bm \xi^l(t_0) = \bm \xi^l_0 = \nabla \bm \pi^{l, *}(\bm \zeta^l_0)$.
Note that in (\ref{alg}), indices $i$ and $j$ are used to differentiate the agents in $\mathcal G_1$ and $\mathcal G_2$. Hereafter, we employ a unifying index $k$ in place of $i$ and $j$, but adopt $l$ to indicate the subnetwork $\mathcal G_l$.
In particular, agent $k$ in network $\mathcal G_l$ means that $k \in \mathcal V_1$ if $l = 1$, and 
$k \in \mathcal V_2$ if $l = 2$.

We aim to improve the robustness of (\ref{alg:asim}) by designing a suitable restarting mechanism.
Compared to \cite{poveda2019inducing, poveda2021robust, ochoa2020robust,ochoa2023momentum}, there are three main challenges: a) (\ref{alg:asim}) is a primal-dual dynamics rather than a purely primal-based dynamics, and the coupling between variables makes the analysis more difficult;
b) the generating functions, beneficial for dealing with constraints, complicate the proof;
c) due to the presence of two subnetworks, it is necessary to design a coordinated mechanism between all agents in $\mathcal G_1$ and $\mathcal G_2$.

We endow agent $k$ in subnetwork $\mathcal G_l$ with a local timer $\tau^l_k$, initialized in the interval $[T_0, T)$. Let the dynamics of $\tau^l_k$ be 
\begin{equation}
\label{timer:dyn}
\tau^l_k \in [T_0, T), ~\dot \tau^l_k =\eta, ~{\rm and}~
\tau^l_k = T, ~\tau^{l, +}_k = T_0,
\end{equation}
where $\eta > 0$ and $T > T_0 > 0$ are manually selected.
Denote $\bm \tau^l = {\rm col}\{\tau^l_1, \dots, \tau^l_{n_l}\}$. 
Then the flow map is designed as
\begin{equation}
\begin{aligned}
\label{alg:sim:flow}
\dot {\bm \xi}^l &= D_{\xi}^l(\bm\tau^l) \big[\nabla \bm\pi^{l, *} (\bm\zeta^l) - \bm\xi^l \big] \\
\dot {\bm \zeta}^l &= {D_{\zeta}^l(\bm\tau^l)} \big[\!-\!\nabla \bm h^l(\bm\xi^l) \!-\! \tilde L^l (\bm\xi^l \!+\! { D_{\zeta}^l}(\bm\tau^l) \dot {\bm\xi}^l) \\
&~~~\!-\! \tilde B^l(\bm\xi^{(3 - l)} \!+\! {D_{\zeta}^l}(\bm\tau^l) \dot {\bm\xi}^{(3 - l)}) \big],\\
\dot{\bm \tau}^l &= \bm \eta^l,
\end{aligned}
\end{equation}
where $D_{\zeta}^l(\bm\tau^l) = {\rm diag}(1/r, 1/r) \otimes {\rm diag}(\bm \tau^l) \otimes I_{p_l}$,
$D_{\xi}^l(\bm\tau^l) = [D_{\zeta}^l(\bm\tau^l)]^{-1}$ and
$\bm \eta^l = \eta \otimes 1_{n_l}$.

Let $\xi^1_k = {\rm col}\{x_k, \lambda_k\}$,
$\zeta^1_k = {\rm col}\{u_k, \gamma_k\}$,
$\xi^2_k = {\rm col}\{y_k, \mu_k\}$,
$\zeta^2_k = {\rm col}\{v_k, \nu_k\}$ and 
$z^l_k = {\rm col}\{\xi^l_k, \zeta^l_k,\tau^l_k\}$.
From the perspective of agent $k$ in $\mathcal G_l$, the flow map is
\begin{equation}
\begin{aligned}
\label{alg:flow}
\tau^l_k&\in [T_0, T) \Rightarrow  
\dot z^l_k = 
\left(
\begin{aligned}
	\dot \xi^l_k\\
	\dot \zeta^l_k\\
	\dot \tau^l_k
\end{aligned}
\right) = F^l_k(z^l_k),
\end{aligned}
\end{equation}
where $F^l_k(z^l_k)$ is defined from the right-hand-side of  (\ref{alg:sim:flow}).
The jump map is designed as
\begin{equation}
\begin{aligned}
\label{alg:jump}
\tau^l_k= T \Rightarrow 
z^{l, +}_k = 
\left(
\begin{aligned}
	\xi_k^{l, +} \\
	\zeta_k^{l, +}  \\
	\tau_k^{l, +}
\end{aligned}
\right) = Q^l_k(z^l_k) := 
\left(
\begin{aligned}
	\xi^l_k \\
	\zeta^l_k \\
	T_0
\end{aligned}
\right).
\end{aligned}
\end{equation}

\subsection{Distributed Coordinated Resets}

In practice, it is important to implement restarting mechanisms in a distributed coordinated manner since each agent is endowed with a local timer \cite{ochoa2020robust, ochoa2023momentum}.
This subsection focuses on designing such a mechanism.

For agent $k$ in subnetwork $\mathcal G_l$, we define a set-valued mapping $R^l_k: \mathbb R \rightarrow \mathbb R$ as
\begin{equation}
\label{coordination:mech}
R^l_k(\tau^l_k) := 
\left\{
\begin{aligned}
&T_0, &&{\rm if}~\tau_k \in [T_0, T_0 + r^l_k), \\
&\{T_0, T\}, &&{\rm if}~\tau_k = T_0 + r^l_k,  \\
&T,    &&{\rm if}~\tau_k \in (T_0 + r^l_k, T),
\end{aligned}
\right.
\end{equation}
with $r^l_k \in \big(0, (T - T_0)/(n_1 + n_2) \big)$.
Let $\mathcal N_k^l$ be the neighbor set of agent $k$ in $\mathcal G_l$, i.e., $(k, \hat k) \in \mathcal E_l, \forall \hat k \in \mathcal N_k^l$,
and $\mathcal N^{(3 - l)}_k$ be the neighbor set of agent $k$ in $\mathcal G_{(3 - l)}$, i.e., $(k, \hat l) \in \mathcal E_{\bowtie}, \forall \hat l \in \mathcal N^{(3 - l)}_k$.
With these preparations, our coordinated resetting rule is designed as follows.

If a local timer $\tau^l_k$ satisfies $\tau^l_k = T$, then
\begin{enumerate}
\item agent $k$ in $\mathcal G_l$ resets its own state $z^l_k$ via (\ref{alg:jump});

\item agent $k$ in $\mathcal G_l$ sends a pulse to its neighbors $\hat k  \in \mathcal N_k^l$, and agent $\hat k$ resets its state $z^l_{\hat k}$ by
\begin{equation*}
\begin{aligned}
	\xi_{\hat k}^{l, +} = \xi_{\hat k}^l, ~
	\zeta_{\hat k}^{l, +}  = \zeta_{\hat k}^l,~
	\tau_{\hat k}^{l, +} \in R^l_{\hat k}(\tau^l_{\hat k});
\end{aligned}
\end{equation*}

\item agent $k$ in $\mathcal G_l$ sends a pulse to its neighbors $\hat l \in \mathcal N^{(3 - l)}_k$, and agent $\hat l$ resets its state $z^{(3 - l)}_{\hat l}$ by 
\begin{equation*}
\begin{aligned}
	\xi_{\hat l}^{(3 - l), +} \!=\! \xi_{\hat l}^{(3 - l)},
	\zeta_{\hat l}^{(3 - l), +}  \!=\! \zeta_{\hat l}^l,
	\tau_{\hat l}^{(3 - l), +} \!\in\! R^{(3 - l)}_{\hat l}
	(\tau^{(3 - l)}_{\hat l});
\end{aligned}
\end{equation*}

\item the remaining agents keep their states.
\end{enumerate}

In contrast to the restarting mechanisms in \cite{ochoa2020robust, ochoa2023momentum}, agents $\hat l \in \mathcal N^{(3 - l)}_k$ are also required to reset their states.
The dynamics (\ref{alg:flow}) and the coordinated resetting mechanism lead to a distributed accelerated algorithm, which is an HDS.
Next, we formalize the system.

Let $\bm z^l = {\rm col}\{z^l_1, \dots, z^l_{n_l}\}$,
$\bm z = {\rm col}\{\bm z^1, \bm z^2\}$,
$\bm \tau = {\rm col}\{\bm \tau^1, \bm \tau^2\}$,
$\bm \xi = {\rm col}\{\bm \xi^1, \bm \xi^2\}$,
$\bm \zeta = {\rm col}\{\bm \zeta^1, \bm \zeta^2\}$,
$\bm \sigma^l = {\rm col}\{\sigma^l_1, \dots, \sigma^l_{n_l}\} \in \mathbb R^{4p_l n_l}$, 
$\bm \sigma = {\rm col}\{\bm \sigma^1, \bm \sigma^2\}$,
$\bm \rho^l = {\rm col}\{\rho^l_1, \dots, \rho^l_{n_l}\} \in \mathbb R^{n_l}$,
$\bm \rho = {\rm col}\{\bm \rho^1, \bm \rho^2\}$, 
$\tilde n = (4p_1 + 1)n_1 + (4p_2 + 1)n_2$,
$\mathbb D_{\xi} = \mathbb D_{\xi}^1 \times \mathbb D_{\xi}^2$ and 
$\mathbb D_{\zeta} = \mathbb D_{\zeta}^1 \times \mathbb D_{\zeta}^2$.
Define a set-valued mapping $G_0: \mathbb R^{\tilde n} \rightrightarrows \mathbb R^{\tilde n}$ as
\begin{equation*}
\begin{aligned}
G_0(&\bm z) = \big\{(\bm \sigma, \bm\rho)\!\in\!\mathbb R^{\tilde n}:
{\rm col}\{\sigma^l_k, \rho^l_k\} = Q^l_k(z_k^l);\\
&\sigma^l_{\hat k} ={\rm col}\{\xi^l_{\hat k}, 
\zeta^l_{\hat k}\},
\rho^l_{\hat k} \in R^l_{\hat k}(\tau^l_{\hat k}), \forall \hat k \in \mathcal N^l_k; \\
&\sigma^{(3-l)}_{\hat l} = {\rm col}\{\xi^{(3-l)}_{\hat l}, \zeta^{(3-l)}_{\hat l}\}, 
\rho^{(3-l)}_{\hat l} \in R^{(3-l)}_k(\tau^l_{\hat l}), \\
&\forall \hat l \in \mathcal N^{(3-l)}_k;
{\rm col}\{\sigma^l_i, \rho^l_i\} = z^l_k, \forall i \in \mathcal V_l \backslash \mathcal N^l_{k}; \\
&{\rm col}\{\sigma^{(3-l)}_j, \rho^{(3-l)}_j\} = z^{(3-l)}_j, \forall j \in \mathcal V_{(3-l)} \backslash \mathcal N^{(3-l)}_{k}\big\}.
\end{aligned}
\end{equation*}

Note that $G_0$ is nonempty only if there is a timer $\tau^l_k$ such that $\tau^l_k = T$, and states of the remaining timers lie in $[T_0, T)$.
Then we design the jump map as
\begin{equation}
\label{dyn:jump}
\bm z^+ \in G_a(\bm z) := \overline G_0(\bm z),
\end{equation}
and the jump set as
\begin{equation}
\label{set:jump}
\mathcal D_a := \big\{\bm z \in \mathbb R^{\tilde n}:
(\bm \xi, \bm \zeta) \in \mathbb D_{\xi} \times \mathbb D_{\zeta},
\max\nolimits_{l, k} \tau^l_k = T\big\},
\end{equation}
where $\overline G_0(\bm z)$ is the outer-semicontinuous hull of $G_0$, i.e., ${\rm graph}(G_a) = {\rm cl}({\rm graph}(G))$.

The flow map $\dot {\bm z} = F_a(\bm z)$ follows from (\ref{alg:flow}), and the flow set is 
$\mathcal C_a := \big\{\bm z \in \mathbb R^{\tilde n}: (\bm \xi, \bm \zeta) \in \mathbb D_{\xi} \times \mathbb D_{\zeta}, \bm \tau \in [T_0, T]^{n_1 + n_2}\big\}$.
Hereafter, we denote the HDS of our distributed accelerated algorithm by $\mathcal H_a = \{\mathcal C_a, F_a, \mathcal D_a, G_a\}$.

\begin{remark}
Compared with the primal-based algorithms in \cite{ochoa2020robust, ochoa2023momentum}, the HDS $\mathcal H_a$ is more powerful
since it is a primal-dual method and can solve constrained optimization/games.
In addition, due to the introduction of generating functions $\bm \psi$ and $\bm \phi$, it covers Nesterov's accelerated scheme \cite{ zeng2022dynamical, ochoa2020robust, ochoa2023momentum} as special cases and handles constraints such as simplices and convex sets efficiently.
\end{remark}

\subsection{Main Results}

This subsection analyzes the performance of $\mathcal H_a$.
Define $\mathcal A = \mathcal  A_{\xi, \zeta} \times \mathcal A_\tau$,
\begin{equation*}
\begin{aligned}
\mathcal A_{\xi, \zeta} =& \big\{(\bm \xi, \bm \zeta) \in \mathbb R^{4p_1n_1 + 4p_2n_2}:
\bm \xi^l = \bm \xi^{l,\star}, 
\bm \xi^{l} = \nabla \bm \pi^{l, \zeta}(\bm \zeta^l)
\big\}, \\
\mathcal A_\tau =& \big\{\bm \tau \in \mathbb R^{n_1 + n_2}:
\bm \tau \in [T_0, T] \cdot 1_{n_1 + n_2} \cup \{T_0, T\}^{n_1 + n_2}\big\},
\end{aligned}
\end{equation*}
where 
$\bm \xi^{1, \star} = {\rm col}\{\bm x^\star, \bm \lambda^\star\}$,
$\bm \xi^{2, \star} = {\rm col}\{\bm y^\star, \bm \mu^\star\}$,
and $(\bm x^\star, \bm \lambda^\star, \bm y^\star, \bm \mu^\star)$ is a saddle point of $S$ in (\ref{lagrangian}).

We first have the following lemma, whose proof is similar to that of Lemma 2 in \cite{ochoa2020robust} or Lemma 1 in \cite{ochoa2023momentum}, and is omitted here due to space limitations.

\begin{lemma}
\label{lem:sol}
Suppose Assumptions \ref{ass:convex} and \ref{ass:genlip} hold. The HDS $\mathcal H_a$ is well-posed in the sense of \cite[Def. 6.29]{goedel2012hybrid}.
Every maximal solution of $\mathcal H_a$ is complete, and there are at most $n_1 + n_2$ jumps at any time interval of length $(T - T_0) / \eta$.
Moreover, $|\bm \tau(t, j)|_{\mathcal A_\tau} = 0$, for all 
$t + j \in \mathcal T(\bm z)$, where 
$$\mathcal T(\bm z)=\{t + j \in {\rm dom}(\bm z): t + j \ge T^*\},$$
and $T^* = (T - T_0)/\eta + (n_1 + n_2)$.
\end{lemma}	

Lemma \ref{lem:sol} indicates that all local timers reach consensus within finite time. As a consequence, the proposed restarting mechanism has the potential to handle faults such as dropped, delayed, or inconsistent reset signals.

Let $\delta > 0$ and $\mathcal A_\delta = \big\{(\bm \xi, \bm \zeta): (\bm \xi, \bm \zeta) \in \mathcal A_{\xi, \zeta} + \delta \mathbb B\big\}$.
We restrict the data of $\mathcal H_a$ by intersecting its $(\bm \xi, \bm \zeta)$ components of flow and jump sets with $\mathcal A_\delta$,
and then, define an HDS as $\mathcal H_{\delta} := \{\mathcal C_{\delta}, F_a, \mathcal D_{\delta}, G_a\}$,
i.e.,
$\mathcal C_{\delta} = \{\bm z \in \mathbb R^{\tilde n}: (\bm \xi, \bm \zeta) \in (\mathbb D_{\xi} \times \mathbb D_{\zeta}) \cap \mathcal A_\delta, \bm \tau \in [T_0, T]^{n_1 + n_2}\}$ and 
$\mathcal D_{\delta} = \{\bm z \in \mathbb R^{\tilde n}: (\bm \xi, \bm \zeta) \in (\mathbb D_{\xi} \times \mathbb D_{\zeta}) \cap \mathcal A_\delta, \max_{k, l} \tau^l_k = T\}$.
It follows from Lemma \ref{lem:sol} that the HDS $\mathcal H_{\delta}$ renders UGAS the compact set $\mathcal A_{\delta, s} = \big((\mathbb D_{\xi} \times \mathbb D_{\zeta}) \cap \mathcal A_\delta \big) \times \mathcal A_\tau$ in finite time.

With the observations, we further restrict the flow and jump sets of $\mathcal H_\delta$ by intersecting them with the set $\mathcal A_{\delta, s}$, and define an HDS as
$\mathcal H_{\delta, s} := \{\mathcal C_{\delta, s}, F_a, \mathcal D_{\delta, s}, G_a\}$, where
$\mathcal C_{\delta, s} = \mathcal C_{\delta} \cap \mathcal A_{\delta, s}$ and 
$\mathcal D_{\delta, s} = \mathcal D_{\delta} \cap \mathcal A_{\delta, s}$.
During the flows, $\bm \tau = \alpha \cdot 1_{n_1 + n_2}$ for $\mathcal H_{\delta, s}$, where $\alpha \in [T_0, T]$.

The following lemma addresses the convergence of $\mathcal H_{\delta, s}$.
\begin{lemma}
\label{lem:HDS:R}
Consider the HDS $\mathcal H_{\delta, s}$ with $\eta = 1$ and $r \ge 2$. 
Under Assumptions \ref{ass:convex} and \ref{ass:genlip},  $\mathcal H_{\delta, s}$ renders UGAS the set $\mathcal A$.
\end{lemma}

\emph{Proof:}
Construct a Lyapunov function candidate as
\begin{equation}
\begin{aligned}
\tilde V(\bm z) \!=&\frac {\alpha^2}{r}\big(S(\bm x, \bm \lambda^\star, \bm y^\star, \bm \mu) \!-\! S(\bm x^\star, \bm\lambda, \bm y, \bm\mu^\star)\big)
\!+\! r D_{\psi_{u}^*}(\bm u, \bm u^\star) \\
&+ r D_{\psi_{\gamma}^*}(\bm \gamma, \bm \gamma^\star) 
+ r D_{\phi_{v}^*}(\bm v, \bm v^\star)
+ r D_{\phi_{\nu}^*}(\bm \nu, \bm \nu^\star).
\end{aligned}
\end{equation}

Note that $\tilde V$ is similar to  $V$ in (\ref{pf:lya0}), where $\alpha$ is used instead of $t$.
Note that the flow map of $\mathcal H_{\delta, s}$ follows from (\ref{alg}). By a similar procedure as in the proof of (\ref{ine:lya}), we can conclude that during flow maps,
\begin{equation}
\label{ine:lya1}
\dot {\tilde V} \le (\frac {2\alpha}{r} - \alpha)\big(S(\bm x, \bm \lambda^\star, \bm y^\star, \bm \mu) - S(\bm x^\star, \bm \lambda, \bm y, \bm \mu^\star)\big).
\end{equation}
Then $\dot {\tilde V} \le 0$ due to $r \ge 2$.
Note that $\dot {\tilde V} = 0$ only if $\bm z \in \mathcal A$.

On the other hand, during the jumps, 
\begin{equation*}
\begin{aligned}
\label{ine:lya2}
\triangle \tilde V(\bm z) =& \tilde V(\bm z^+) - \tilde V(\bm z) 
= - \big(\frac {T^2}{r} - \frac {T_{\min}^2}{r}\big)\\ 
&\cdot \big(S(\bm x, \bm \lambda^\star, \bm y^\star, \bm \mu) 
- S(\bm x^\star, \bm\lambda, \bm y, \bm\mu^\star)\big) \le 0.
\end{aligned}
\end{equation*}

In summary, $\tilde V$ is strictly decreasing during flows and does not increase during jumps. By the hybrid invariance principle \cite[Theorem 8.8]{goedel2012hybrid}, 
$\mathcal H_{\delta, s}$ renders UGAS the set $\mathcal A$ and the proof is complete.
$\hfill\square$

The next result establishes the convergence of $\mathcal H_a$.

\begin{theorem}
\label{thm:convergence}
Consider the HDS $\mathcal H_a$ with $\eta = 1$ and $r \ge 2$. Under Assumptions \ref{ass:convex} and \ref{ass:genlip},
\begin{enumerate}[i)]
\item the set $\mathcal A$ is R-UGAS;

\item during flows, it holds that 
\begin{equation*}
	\begin{aligned}
		S \big(\bm x(t, j),\bm \lambda^*, \bm y^*, \bm \mu(t, j) \big) - S \big(\bm x^*, \bm \lambda(t, j), \bm y(t, j), \bm \mu^*\big)&\\
		\le \frac {c_j}{\big[\tau^l_k(t, j)\big]^2},~
		\forall t + j \in \mathcal T(\bm z)&,
	\end{aligned}
\end{equation*}
where $\{c_j\}_{j = 0}^\infty$ is a monotonically decreasing sequence of positive numbers and approaches $0$.
\end{enumerate}
\end{theorem}

\emph{Proof:}
i) Since $\mathcal A$ is UGAS for $\mathcal H_{\delta, s}$, and $\mathcal A_{\delta, s}$ is UGAS in finite time for $\mathcal H_{\delta}$, it follows from the hybrid reduction principle \cite[Corollary 7.24]{goedel2012hybrid} that
$\mathcal A$ is UGAS for $\mathcal H_{\delta}$.
Due to the arbitrariness of $\delta$, $\mathcal A$ is also UGAS for $\mathcal H_a$.
In light of \cite[Theorem 7.21]{goedel2012hybrid}, $\mathcal A$ is R-UGAS for $\mathcal H_a$ by the well-posedness of $\mathcal H_a$. Thus, part i) is proved.

ii) Recall from Lemma \ref{lem:sol} that $|\bm \tau(t, j)|_{\mathcal A_\tau} = 0$ for all $t + j \in \mathcal T(\bm z)$.
For each $(t, j), (s, j)\in \mathcal T(\bm z)$ with $t > s$, 
$\tilde V(\bm z(t, j)) \le \tilde V(\bm z(s, j))$ since $\dot{\tilde V} \le 0$.
As a consequence, 
\begin{equation*}
\begin{aligned}
&S \big(\bm x(t, j),\bm \lambda^*, \bm y^*, \bm \mu(t, j)\big)
- S\big(\bm x^*, \bm \lambda(t, j),\bm y(t, j), \bm \mu^*\big)\\
&\le \frac{r}{\big[\tau^l_k(t, j)\big]^2} \tilde V(\bm z(t, j))
\le \frac{r}{\big[\tau^l_k(t, j)\big]^2} \tilde V(\bm z(s_j, j))
= \frac{c_j}{\alpha^2},
\end{aligned}
\end{equation*}
where $c_j = r \tilde V(\bm z(s_j, j))$ and $s_j = \min\{t \ge 0: (t, j) \in \mathcal T(\bm z)\}$.
Since $\tilde V$ is nonincreasing and  converges to $0$, $\{c_j\}_{j = 0}^\infty$ is a monotonically decreasing sequence, and approaches $0$.
The proof is complete.
$\hfill\square$

Theorem \ref{thm:convergence} implies that the proposed $\mathcal H_a$ is UGAS with structural robustness, and also enjoys a rate of $\mathcal O(1/t^2)$. 

\begin{example}
Consider solving the problem in Example $1$ by the HDS $\mathcal H_a$.
Fig. \ref{fig:Ex2} shows the trajectories of $U(\bm x, \bm y)$ and indicates that $\mathcal H_a$ is convergent under disturbances.
\begin{figure}[htp]
\centering
\includegraphics[scale=0.4]{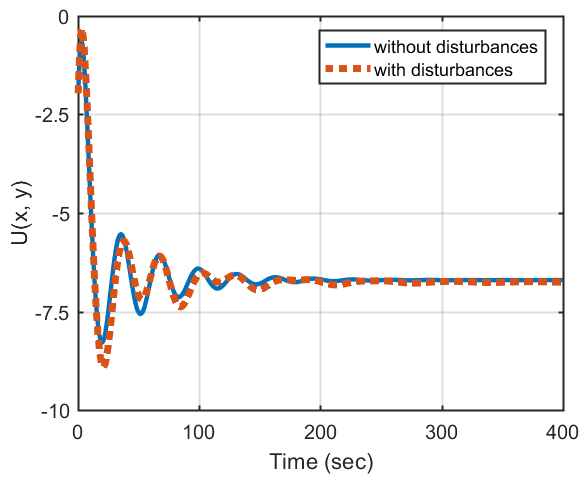}
\caption{Trajectories of $U(\bm x, \bm y)$ under the HDS $\mathcal H_a$.
}
\label{fig:Ex2}
\end{figure}
\end{example}

\section{CONCLUDING REMARKS}
In this paper, we studied the problem of finding an NE in a two-subnetwork bilinear zero-sum game. We first proposed a primal-dual accelerated algorithm, which, however, failed to converge under disturbances. To address this, we introduced a restarting mechanism and developed a distributed algorithm that achieves exact NE with structural robustness. As a future direction, we aim to extend our results to general-sum and $N$-player monotone games to enhance their applicability.

\bibliographystyle{IEEEtran}
\bibliography{references,IDS_Publications_03112025}

\end{document}